\documentclass[12pt, reqno]{amsart}
\usepackage{graphicx}
\usepackage{color}
\usepackage[english]{babel}
\usepackage[latin1]{inputenc}
\usepackage{amsmath}
\usepackage{amssymb}
\newtheorem{teo}{Theorem}[section]

\newtheorem{lema}[teo]{Lemma}
\newtheorem{prop}[teo]{Proposition}
\theoremstyle{definition}

\theoremstyle{remark}

\numberwithin{equation}{section}
\newcommand{\To}{\longrightarrow}

\newcommand{\N}{\mathbb{N}}

\newcommand{\C}{\mathbb{C}}
\newcommand{\D}{\mathbf{D}}
\newcommand{\dem}{\noindent {\textsf{Proof.} }}
\newcommand{\Hinf}{H^{\infty}(B_E)}
\renewcommand{\phi}{\varphi}

\newcommand{\ra}{\rangle}
\newcommand{\la}{\langle}
\author[Miralles]{Alejandro Miralles \\ 
\textit{IMAC \lowercase{and} D\lowercase{epartament de} M\lowercase{atemàtiques} \\
U\lowercase{niversitat} J\lowercase{aume} I \lowercase{de} C\lowercase{astelló} (UJI), C\lowercase{astelló} E-12071, Spain.\\
\lowercase{\emph{e}.mail: mirallea@uji.es,} \ T\lowercase{el:} +34 964728411.}}

\thanks{The author is supported by Project MTM 2011-22457 (Ministerio de Economía y Competitividad. Spain) and Project P1-1B2014-35 (Universitat Jaume I. Spain).}


\begin{document}

\title[Interpolating sequences for $H^{\infty}(B_H)$]{Interpolating sequences for $H^{\infty}(B_H)$}

\begin{abstract}
We prove that under the extended Carleson's condition, a sequence $(x_n) \subset B_H$ is linear interpolating for $H^{\infty}(B_H)$ for an infinite dimensional Hilbert space H. In particular, we construct the interpolating functions for each sequence and find a bound for the constant of interpolation. 
\end{abstract}

\maketitle

\section{Introduction}

Let $A$ be a space of bounded functions defined on $X$. A sequence $(x_n)$ in X
is called interpolating for $A$ if for
 any sequence  $(a_n)\in \ell_\infty,$ there exists $f\in A $ such that
 $f(x_n)=a_n \hbox{ for all} \;n\in \mathbb{N}.$ We consider the linear and continuous $R: A \to \ell_{\infty}$ defined by $R(f)=(f(x_n))$. The sequence $(x_n)$ is interpolating for $A$ if and only if there exists a map $T: \ell_\infty \to A$ such that $R \circ T = Id_{\ell{\infty}}$. If $T$ is linear, the sequence $(x_n)$ is said to be linear interpolating for $A$. For any $\alpha=(\alpha_j) \in \ell^{\infty},$ let  $M_{\alpha}=\inf
\left\{ \|f\|_{\infty} : f(x_j)= \alpha_j, \; j \in
\mathbb{N}, f \in A \right\}.$ The constant of interpolation for $(x_n)$ is defined by $M= \sup \left\{ M_\alpha : \alpha \in \ell^{\infty}, \| \alpha \|_{\infty} \leq 1 \right\}.$
%Any real number not smaller than $M$ is a \emph{constant of interpolation} for
%$(x_n).$

It is a classical result in function theory that a sequence $(z_n)$
in the open unit disc $\mathbf{D}\subset \mathbb{C}$ is
interpolating for $H^\infty ,$ the space of analytic bounded functions on $\D$, if and only if Carleson's condition
holds, i. e.:
\begin{eqnarray} \label{carleson}
\mbox{ There is } \delta > 0 \mbox{ such that } \  \prod_{k\neq
j} \rho(z_k,z_j) \geq
\delta \ \mbox{ for any }j \in \mathbb{N},
\end{eqnarray}

\noindent where $\rho(z_k,z_j)$ denotes the pseudohyperbolic distance for points $z_k,z_j \in \D$, given by $$\rho(z_k,z_j)=\left|\frac{z_{k}-z_{j}}{1-\overline{z_{k}}z_{j}}\right|.$$

Recall the Schwarz-Pick lemma: $\rho(f(z),f(w)) \leq \rho(z,w)$ for any $z, w \in \D$ and $f \in H^\infty$, $\|f\| \leq~1$. If $\psi$ is an automorphism of $\D$, then  $\rho(\psi(z),\psi(w))=\rho(z,w)$.

 If we deal with complex Banach spaces $E$, a function $f:B_E \to \mathbb{C}$ is said to be analytic if it is Fr{\'e}chet differentiable. Denote
by $\Hinf$  the space $\{f:B_E \to \mathbb{C}: f\hbox{ is analytic
and bounded } \}$, which becomes a uniform Banach algebra when endowed with
the sup-norm $\|f\|=\sup\{|f(x)|:x\in B_E\}$ and it is, obviously,
the analogue of the space $H^\infty$ for an arbitrary Banach space.

Sufficient conditions for a sequence to be interpolating for $H^\infty(B_E)$ where given by the authors in \cite{GM}. Bearing in mind the Davie-Gamelin extension of $f \in H^{\infty}(B_E)$ to $\tilde{f} \in H^{\infty}(B_{E^{**}})$, the authors proved that a sufficient condition for a sequence $(x_n) \subset B_{E^{**}}$ to be linear interpolating for $\Hinf$ is that the sequence of norms $(\|x_n\|)$ is interpolating for $H^\infty$. Examples of sequences which satisfy this condition are, for instance, those which grow exponentially to the unit sphere, which we call the Hayman-Newman condition: $1-\| x_{k+1}\| < c(1-\| x_{k}\|)$ for some $0 < c < 1$ for any $k \in \N$. Interpolating sequences on $H^{\infty}(B_E)$ has been very useful to study the spectra of composition operators on spaces of analytic functions (see \cite{GGL}, \cite{GLR} and \cite{GM2}).

The notion of pseudohyperbolic distance can be carried over to $H^{\infty}(B_H)$ by considering for any $x, y \in B_H$,
\begin{eqnarray} \label{ref pseud}
\rho_{H}(x,y)= \sup \{ \rho(f(x),f(y)) : f \in H^{\infty}(B_H), \|f\| \leq 1 \},
\end{eqnarray}

\noindent where $\rho (z,w)$ is the pseudohyperbolic distance in $\D.$  B. Berndtsson \cite{Ber85}  showed that a sequence
$(x_n)$ in the open unit Euclidean ball $B_n$ of $\mathbb{C}^n$ is
interpolating for  $H^\infty(B_n)$
 if the
following extended Carleson's condition holds: 
\begin{eqnarray} \label{carleson1} 
\mbox{ There is }
\delta > 0 \mbox{ such that }\qquad \prod_{k \neq j} \rho_H (
x_j,x_k) \geq \delta \qquad \forall  j \in \mathbb{N}.
\end{eqnarray} 

%where
%$\rho_H$ denotes the pseudohyperbolic distance in $B_n.$

As P. Galindo, T. Gamelin and M. Lindstr{\"o}m  pointed out in \cite{GGL}, the result given by Berndtsson can be extended to the case of an infinite dimensional complex Hilbert space $H$ by interpolating on finite subsets of the sequence with uniform bounds and applying a normal families argument. 

The aim of this paper is to adapt the proof given by Berndtsson to the infinite dimensional case and prove that under the extended Carleson's condition \ref{ref pseud}, a sequence $(x_n) \subset B_H$ is linear interpolating. In particular, we will construct the interpolating functions for each sequence and will find a bound for the constant of interpolation. 

For our purpose, we will study the automorphisms on $B_H$ and will adapt some results given by B. Berndtsson to the infinite dimensional case. 

\section{Background}

\textbf{Automorphisms on $B_H$.} Recall that the set of automorphisms on $\D$ is denoted by $Aut(\D)$. It is well-known that this set is generated by rotations and Möbius transformations $m_a: \D \To \D$ given by
\begin{eqnarray} \label{automorph2}
m_a(z)=\frac{a-z}{1-\bar{a}z} \mbox{ \ for any $a \in \D$.}
\end{eqnarray}

The analogues of M\"{o}bius transformations on $H$ are $\phi_a: B_{H} \To B_{H}\;, a \in B_{H},$ defined according to
\begin{eqnarray} \label{automorph}
\phi_{a}(x)= (s_a Q_a+P_a)(m_a(x))
\end{eqnarray}

\noindent where $s_a=\sqrt{1-\|a\|^2},$  $m_a: B_{H} \To B_{H}$ is  the analytic map
\begin{eqnarray} \label{mobius hilbert}
m_a(x)=\frac{a-x}{1- \la x,a \ra},
\end{eqnarray}

\noindent  $P_a: H \To H$ is the orthogonal projection along the one-dimensional subspace spanned by $a$, that is,
$$P_a(x)=\frac{\langle x, a \rangle}{\langle a,a \rangle} a$$

\noindent and $Q_a: H \To H,$  is the orthogonal complement, $Q_a=Id-P_a$. Recall that $P_a$ and $Q_a$ are self-adjoint operators since they are orthogonal projections, so $\langle P_a(x),y \rangle=\langle x,P_a(y) \rangle$ and $\langle Q_a(x),y \rangle=\langle x,Q_a(y) \rangle$ for any $x, y \in H$.

The automorphisms of the unit ball $B_{H}$ turn to be compositions of such analogous
M\"{o}bius transformations with  unitary transformations $U$ of $H$, that is,  self-maps of $H$
 satisfying $\langle U(x) , U(y) \rangle = \langle x, y \rangle$ for all $x, y \in H$. \medskip

\textbf{Remarks on the pseudohyperbolic distance.} It is clear by the definition that for any $x,y \in B_H$,
\begin{eqnarray} \label{metrica pseudohip contractiva}
\rho (f(x),f(y)) \leq \rho_{H}(x,y)~~ \text{ for any } f \in H^{\infty}(B_{H}), \ \|f\|_{\infty} \leq 1,
\end{eqnarray} 

It is also well-known that 
\begin{eqnarray} \label{metrica hilb}
\rho_H(x,y) = \|\phi_{y}(x)\|
\end{eqnarray}

\noindent so making some calculations we obtain
\begin{eqnarray} \label{pseudometric Hilbert}
\rho(x,y)^2= 1 - \frac{(1-\|x\|^2)(1-\|y\|^2)}{|1- \la x,y \ra|^2}.
\end{eqnarray}

\noindent All these results and further information on the automorphisms of $B_{H}$ and the pseudohyperbolic distance can be found in \cite{GR}. \medskip

\section{Results}
First, we recall Proposition 2.1 in \cite{GLM09}:
\begin{prop} \label{prop 21}
Let $E$ be a complex Banach space and $(x_n) \subset B_E$. If there exists $M >0$ and a sequence of functions $(F_j)\subset H^{\infty}(B_E)$ satisfying $F_j(x_n)=\delta_{j,n}$ for any $j \in \N$ and $\sum_j |F_j(x)|\leq M$ for
all $x\in B_E$, then $(x_n)$ is linear interpolating for $H^{\infty}(B_E)$.
\end{prop}

We will call $(F_n)$ a sequence of Beurling functions for $(x_n)$. Under conditions of Proposition \ref{prop 21}, we have that $T:\ell_\infty \to H^{\infty}(B_E)$ defined by
$T((\alpha_n))=\sum_n \alpha_n F_n$ is a well-defined, linear operator such
that $\|T\|\leq M$ and $T((\alpha_n))(x_k)=\alpha_k$ for any
$k\in \mathbb{N},$ so $(x_n)$ is linear interpolating. In particular, the constant $M$ is an upper bound for the constant of interpolation. \medskip

The following calculations are straightforward and can be found in \cite{GM}.

\begin{lema} \label{lemametps}
We have the following statements: \medskip
%\begin{itemize}
%\begin{eqnarray} \label{lemametps1}
%a) \qquad \rho(a,c) \geq \rho (b,c) \quad \mbox{ for real numbers } 0 \leq a \leq b
%\leq c < 1.
%\end{eqnarray}
\begin{eqnarray} \label{lemametps2}
\qquad 1-x \leq - \log x \quad \mbox{ for } 0 < x \leq 1. \hspace{4.15cm}
\end{eqnarray}
\begin{eqnarray} \label{lemametps3}
\qquad \Re e \ \left[ \frac{1+ \alpha z}{1- \alpha z} \right] =
\frac{1-|\alpha|^2 |z|^2}{|1-\alpha z|^2} \quad \mbox{ for any } \alpha \in \overline{\D}, z \in \D.
\end{eqnarray}
\end{lema} \medskip

%First, we recall Lemma 3 from \cite{Ber85}:
%\begin{lema}
%Consider a sequence $(S_j)$ of sets of $\N$ satisfying the property that if $k \notin S_j$, then $j \in S_k$. Then,
%\begin{eqnarray} 
%\sum_{j=1}^{\infty} c_j h \left( \sum_{k \in S_j} c_k \right) \leq 2 e.
%\end{eqnarray}
%\end{lema}

%We will also need the following result:
%\begin{lema} \label{lema ber hilbert 2}
%Let $h(t)= \min \left\{ 1, 1/t^2 \right\}$. Then, the function $x^2 \exp{\left( -xt/8 \right)}$ is bounded by $256 h(t)/e^2$ \ for \ $0 \leq x \leq 1$ and $t > 0$.
%\end{lema}

The following three lemmas are just calculus:
\begin{lema} \label{lema ber hilbert 2}
The function $u^2 \exp{\left( -u t/8 \right)}$ is upper bounded by $\min\{ 1,\frac{256}{e^2 t^2} \}$ \ for \ $0 \leq u \leq 1$ and $t > 0$.
\end{lema}

%The following three lemmas are just calculations:
%\begin{lema} \label{lema ber hilbert 2}
%The function $x^2 \exp{\left( -xt/8 \right)}$ is upper bounded by $\min\{ 1,\frac{128}{t^2} \}$ \ for \ $0 \leq x \leq 1$ and $t > 0$.
%\end{lema}
%
%
%\dem Bearing in mind the Taylor series of $\exp{x}$, we have 
%$$x^2 \exp{\left( -xt/8 \right)}=\frac{x^2}{\exp{\left( xt/8 \right)}}=\frac{x^2}{1+xt/8+x^2t^2/128+ \cdots} \leq \min\{ 1, 128/t^2 \}.$$
%
%\begin{lema} \label{lema ber hilbert 3}
%We have
%$$\int_{0}^{\infty} \min\{1,\frac{128}{t^2}\}dt= 16 \sqrt{2}$$
%\end{lema}

%Next lemma is a particular case of Lemma 3 in \cite{BCL}:
\begin{lema} \label{func dec}
Let $0<c_k <1$ for any $k \in \N$ and suppose that $h(t)$ is a non-increasing function on $(0,\infty)$. Then,
$$\sum_{j=1}^n c_j h\left( \sum_{k \geq j} c_j \right) \leq \int_{0}^{\infty} h(t) dt.$$
%Suppose that $S_j=\{ k \in \N : \|x_k\| \geq \|x_j\| \}$. Then
%$$\sum_{k=1}^{\infty} c_k \exp(-\sum_{k \in S_j} c_j) \leq 2e.$$
\end{lema}

The following result will be needed to simplify the proof of Theorem \ref{teo carl hilbert}. Maybe it is folklore but we prove it for the sake of completeness:
\begin{lema} \label{reordenar}
Let $(a_n)\subset [0,1)$ such that $\lim_n a_n=1$. Then, $(a_n)$ can be reordered into a non-decreasing sequence $(b_n)$ such that $\lim_{n} b_n=1$.
\end{lema}

\dem Consider $m \in \N$. Since $\lim_n a_n=1$, there exists $n_m \in \N$ such that $a_n \geq 1-\frac{1}{m+1}$ for any $n \geq n_m$, so the set $B_m=\{a_n : 1-\frac{1}{m} \leq a_n < 1-\frac{1}{m+1} \}$ is finite. It is clear that there exists $m_0 \in \N$ such that for any $m \geq m_0$, the set $B_m$ is non-avoid and $B_{m_0-1}$ is avoid (where $B_0=\varnothing$). Since $B_{m_0}$ is finite, we order its elements: $b_1,b_2, \cdots, b_{r_0}$. Now, we order the elements of $B_{m_0+1}$: $b_{r_0+1}, b_{r_0+2}, \cdots, b_{r_1}$ and so on. It is clear that we reorder the whole sequence $(a_n)$: Given $a_n$ for some $n \in \N$, there exists $m_n$ such that $1-\frac{1}{m_n} \leq a_n < 1-\frac{1}{m_n+1}$, so it belongs to $B_{m_n}$. \qed \bigskip

Now we provide a lemma which includes some calculations related to the automorphisms $\phi_{x}$. 

\begin{lema} \label{lemacarlesonhilbert}
Let $x,y \in B_H$ and $\phi_{-y}:H \To H$ the corresponding automorphism defined as in \ref{automorph}. Then, we have that
$$\la \phi_{y}(x),\phi_{y}(z)\ra=1- \frac{(1-\la x,z\ra)(1-\la y,y\ra)}{(1-\la x,y\ra)(1-\la y,z\ra)}.$$
\end{lema}

\dem Since for any $x \in B_H$ we have $\phi_{y}(x)= \left( s_y Q_{y}+P_{y} \right)(m_{y}(x) ),$ and bearing in mind that $P$ and $Q$ are orthogonal, we obtain that
$$\langle \phi_{y}(x),\phi_{y}(z)\rangle = \langle \left( s_y Q_{y}+P_{y} \right)(m_{y}(x)), \left( s_y Q_{y}+P_{y} \right)(m_{y}(z)) \rangle=$$
$$s_y^2 \langle Q_{y}(m_{y}(x)), Q_{y}(m_{y}(z)) \rangle + \langle P_{y}(m_{y}(x)),P_{y}(m_{y}(z)) \rangle=$$
$$\frac{(1-\|y\|^2) \langle Q_{y}(y-x), Q_{y}(y-z) \rangle+ \langle P_{y}(y-x),P_{y}(y-z) \rangle}{(1- \langle x,y \rangle)\overline{(1- \langle z,y \rangle)}}$$

\noindent by (\ref{mobius hilbert}) just making some calculations. Since we have that $P_a + Q_a =Id_H$ for any $a \in H$, we obtain that 
$$\langle \phi_{y}(x),\phi_{y}(z) \rangle = \frac{\langle y-x ,y-z \rangle-\|y\|^2 \langle Q_{y}(y-x), Q_{y}(y-z) \rangle}{(1- \langle x,y \rangle)(1- \langle y,z \rangle)}.$$

The complement of the orthogonal projection is given by
$$Q_{y}(x)=x-\frac{\langle x,y \rangle}{ \langle y,y \rangle}y,$$

\noindent hence $Q_{y}(y-x)=-Q_{y}(x)$ and $Q_{y}(y-z)=-Q_{y}(z)$. Moreover,

\begin{eqnarray*}
\langle -Q_{y}(x), -Q_{y}(z) \rangle= \langle Q_{y}(x), Q_{y}(z) \rangle=  \langle x-\frac{\la x,y \rangle}{ \langle y,y \rangle}y,z-\frac{\langle z,y \rangle}{ \langle y,y \rangle}y \rangle=\\
\langle x,z\rangle-\frac{1}{\|y\|^2}\langle x,y \rangle \langle y,z \rangle-\frac{1}{\|y\|^2}\langle x,y \rangle \langle y,z \rangle+\frac{1}{\|y\|^2} \langle x,y \rangle \langle y,z \rangle=\\
\langle x,z \rangle -\frac{1}{\|y\|^2} \langle x,y \rangle  \langle y,z \rangle= \frac{\langle x,z \rangle  \langle y,y \rangle-\langle x,y \rangle \langle y,z \rangle}{\|y\|^2}.
\end{eqnarray*}

\begin{eqnarray*}
\mbox{So, } \langle \phi_{y}(x),\phi_{y}(z)\rangle= \frac{\langle x-y,z-y \rangle-\|y\|^2\frac{\langle x,z \rangle \langle y,y \rangle- \langle x,y\rangle \langle y,z\rangle}{\|y\|^2}}{(1- \langle x,y\rangle)(1- \langle y,z\rangle)}= \\
\frac{ \langle x-y,z-y\rangle- \langle x,z\rangle \langle y,y\rangle+ \langle x,y\rangle \langle y,z\rangle}{(1- \langle x,y\rangle)(1-\langle y,z\rangle)}.
\end{eqnarray*}

Adding and subtracting $1$ and arranging terms, we obtain that the numerator equals to $(1-\la x,y \ra)(1-\la y,z \ra)-(1-\la x,z \ra)(1-\la y,y \ra).$

Therefore, dividing by the denominator and making calculations, we obtain
$$\la \phi_{y}(x),\phi_{y}(z) \ra=1-\frac{(1-\la x,z \ra)(1-\la y,y \ra)}{(1- \la x,y \ra)(1- \la y,z \ra)},$$
%\frac{(1-<x,y\ra)(1-<y,z>)-(1-<x,z>)(1-<y,y>)}{(1-<x,y>)(1-<y,z>)}=$$
\noindent and the lemma is proved. \qed \bigskip

Considering $z=x$, we obtain that formulas \ref{metrica hilb} and \ref{pseudometric Hilbert} are the same expression for the pseudohyperbolic distance for $x,y \in B_H$. \medskip

We will also need some technical lemmas. For the first one, we will need Proposition 5.1.2 in \cite{R4}, which is stated as follows,

\begin{lema} \label{lema rudin}
Let $a, b, c$ be points in the unit ball of a finite dimensional Hilbert space. Then,
$$|1- \la a,b\ra| \leq  (\sqrt{|1-\la a,c\ra|} + \sqrt{|1-\la b,c\ra|})^2$$
%and
%$$1-|<a,b>| \leq 2(1-|<a,c>| + 1-|<b,c>|).$$
\end{lema}

Then, we obtain the following lemma, %which is an extension of Lemma 5 in \cite{Ber85} to any complex Hilbert space,

\begin{lema} \label{desigualdad hilb}
Let $H$ be an infinite dimensional complex Hilbert space and $x_1,x_2,x_3 \in B_H$. Then,
$$|1-\la x_1,x_2\ra| \leq 2 (|1-\la x_1,x_3\ra| + |1-\la x_2,x_3\ra|)$$

\noindent and
$$1-|\la x_1,x_2\ra| \leq 2 (1-|\la x_1,x_3\ra| + 1-|\la x_2,x_3\ra|).$$
\end{lema}

\dem Let $x_1,x_2,x_3 \in \overline{B}_H$ and set the space $H_1 =$ span$\{ x_1,x_2,x_3 \}$. We have that $H_1$ is itself a Hilbert space and we can consider an orthonormal basis $\left\{ e_1, e_2, e_3 \right\}$ of $H_1$. Consider for $j=1,2,3$ the vectors $y_j=(y_j^1, y_j^2,y_j^3)$ given by the components of $x_j$ in that basis. It is clear that these vectors are in the unit Euclidean ball of $\C^3$ and $\la x_j,x_k\ra=\la y_j,y_k\ra$, so we apply Lemma \ref{lema rudin} to deduce
$$|1- \la x_1,x_2\ra| \leq  (\sqrt{|1-\la x_1,x_3\ra|} + \sqrt{|1-\la x_2,x_3\ra|})^2 =$$
$$|1-\la x_1,x_3\ra|+|1-\la x_2,x_3\ra|+ 2 \sqrt{|1-\la x_1,x_3\ra|} \sqrt{|1-\la x_2,x_3\ra|}.$$

%Let $a,b,c \in \overline{B}_H$ and set the tridimensional space $H_1 = <a,b,c>$. We have that $H_1$ is itself a Hilbert space and it is contained in $H$ isometrically, so we apply Lemma \ref{lema rudin} to deduce
%$$|1- <a,b\ra| \leq  (\sqrt{|1-<a,c\ra|} + \sqrt{|1-<b,c\ra|})^2 =$$
%$$|1-<a,c\ra|+|1-<b,c\ra|+ 2 \sqrt{|1-<a,c\ra|} \sqrt{|1-<b,c\ra|}.$$

\noindent By the arithmetic-geometric means inequality, we have $|1-\la x_1,x_2\ra~|~\leq$
$$|1-\la x_1,x_3\ra|+|1-\la x_2,x_3\ra| +2 \frac{|1-\la x_1,x_3\ra|+|1-\la x_2,x_3\ra|}{2} =$$
$$ 2 (|1-\la x_1,x_3\ra|+|1-\la x_2,x_3\ra|).$$

\noindent To prove the other result, notice that
$$1-|\la x_j,x_k\ra| = \min_{\theta \in [0, 2 \pi)} |1-e^{i \theta} \la x_j,x_k\ra|.$$

We have that $1-|\la x_1,x_3\ra| = |1-e^{i \alpha} \la x_1,x_3\ra|$ and
$1-|\la x_2,x_3\ra| = |1-e^{i \beta} \la x_2,x_3\ra|$ for some $\alpha, \beta \in [0, 2\pi)$. Then, applying the inequality above, we have that
\begin{eqnarray*}
1-|\la x_1,x_2\ra| = 1-| \la e^{i \alpha} x_1, e^{i \beta} x_2\ra| \leq 
|1-\la e^{i \alpha} x_1, e^{i \beta} x_2\ra| \leq \\
2( |1-e^{i \alpha} \la x_1,x_3\ra| + |1-e^{i \beta} \la x_2,x_3\ra|) = \\
 2(1-|\la x_1,x_3\ra| + 1-|\la x_2,x_3\ra|). \qed
\end{eqnarray*}  \smallskip

Then, we can prove the following lemma, which extends Lemma 6 in \cite{Ber85} to the infinite dimensional case. We give the proof for the sake of completeness. 

\begin{lema} \label{lema ber hilbert}
Let $H$ be a Hilbert space and $x_k, x_j \in B_H$. If $\|x_k\| \geq \|x_j\|$, then
\begin{eqnarray} \label{des ber hilbert}
\frac{1-|\la x_k,x\ra|^2}{1-|\la x_k,x_j\ra|^2} \geq \frac{1}{8} \frac{1-\|x_k\|^2}{1-|\la x_j,x\ra|^2} \qquad \mbox{for any } x \in B_H.
\end{eqnarray}
\end{lema}

\dem By Lemma \ref{desigualdad hilb}, $1-|\la x_k,x_j \ra| \leq 2(1-|\la x_k,x \ra|+1-|\la x,x_j \ra|)$, and we consider two cases depending on $x \in B_H$. If $1-|\la x_k,x \ra| \geq 1-|\la x_j,x \ra|$, then $1-|\la x_k,x_j \ra|^2 \leq 8 (1-|\la x_k,x \ra|)$ so, bearing in mind that $\|x_k\| \geq \|x_j\|$, 
$$\frac{1-|\la x_k,x \ra|^2}{1-|\la x_k,x_j \ra|^2} \geq \frac{1}{8} \frac{1-|\la x_k,x \ra|^2}{1-|\la x_k,x \ra|} \geq \frac{1}{8} \geq   \frac{1}{8} \frac{1-\|x_j\|^2}{1-|\la x_j,x \ra|^2} \geq   \frac{1}{8} \frac{1-\|x_k\|^2}{1-|\la x_j,x \ra|^2}.$$

\noindent If $1-|\la x_k,x \ra| \leq 1-|\la x_j,x \ra|$, then $1-|\la x_k,x_j \ra|^2 \leq 8 (1-|\la x_j,x \ra|)$ so,
$$\frac{1-|\la x_k,x \ra|^2}{1-|\la x_k,x_j \ra|^2} \geq \frac{1}{8} \frac{1-|\la x_k,x \ra|^2}{1-|\la x_k,x \ra|} \geq  \frac{1}{8} \frac{1-\|x_k\|^2}{1-|\la x_k,x \ra|} \geq \frac{1}{8} \frac{1-\|x_k\|^2}{1-|\la x_j,x \ra|^2},$$

so we are done. \qed \bigskip

We will also need the following lemma,
\begin{lema} \label{lema4Hilbert}
Let $\{x_n\} \subset B_{H}$ and $\delta > 0$  satisfying
\begin{eqnarray} \label{carleson condition hilbert}
\prod_{k \neq j} \rho(x_k,x_j) \geq \delta \mbox{ for all } j \in \N.
\end{eqnarray}

\noindent Then, we have that
\begin{eqnarray} \label{eqlemaHilbert}
\sum_{k \neq j}^{\infty} (1-\|x_{k}\|^2) \leq 2 \log
\left( \frac{1}{\delta} \right) \frac{1+\|x_{j}\|}{1-\|x_{j}\|} \qquad \mbox{ for any } 
 j \in \mathbb{N},
\end{eqnarray}

\noindent and for any $j \in \N$,
\begin{eqnarray} \label{eqlemaHilbert2}
\sum_{k =1}^{\infty} (1-\|x_{k}\|^2) \leq \left( 1+2 \log
\frac{1}{\delta} \right) \frac{1+\|x_{j}\|}{1-\|x_{j}\|}.
\end{eqnarray}
\end{lema}

\dem  Taking squares and logarithms in \ref{carleson condition hilbert} we obtain
$$- \sum_{k \neq j}^{\infty} \log \rho(x_k,x_j)^2 \leq -2 \log \delta=2 \log{\frac{1}{\delta}}.$$

By (\ref{lemametps2}), we have that $1- \rho(x_k,x_j)^2 \leq - \log \rho(x_k,x_j)^2 \mbox{ for any } k \neq j$, so bearing in mind (\ref{pseudometric Hilbert}), we obtain
$$\sum_{k \neq j}^{\infty} \frac{(1-\|x_{k}\|^2)(1-\|x_{j}\|^2|)}{|1-\la x_{k},x_{j}\ra|^2} \leq
 2 \log{\frac{1}{\delta}}.$$

In consequence,
$$\sum_{k \neq j}^{\infty} (1 - \|x_{k}\|^2) = \sum_{k \neq j}^{\infty} \frac{(1 - \|x_{k}\|^2)(1-\|x_{j})\|^2)}{|1-\la x_{k},x_j\ra|^2}
\frac{|1-\la x_{k},x_j\ra|^2}{1-\|x_j\|^2} \leq$$
$$2 \left(\log{\frac{1}{\delta}}\right) \frac{(1+ \|x_{j}\|)^2}{1-
\|x_{j}\|^2}=2 \left( \log \frac{1}{\delta} \right) \frac{1+ \|x_{j}\|}{1-
\|x_{j}\|}$$
\noindent and the lemma is proved. \hfill $\square$ \bigskip

Now we are ready to prove the result for complex Hilbert spaces. In addition, we will provide an upper estimate for the constant of interpolation depending only on $\delta$. %During the proof, the generalized Carleson condition refers to condition \ref{carleson condition hilbert}.

\begin{teo} \label{teo carl hilbert}
Let $H$ be a Hilbert space and $(x_n)$ a sequence in $B_H$. Suppose that there exists $\delta > 0$ such that $(x_n)$ satisfies the generalized Carleson condition (\ref{carleson1}) for $\delta$. Then, there exists a sequence of Beurling functions $(F_n)$ for $(x_n)$. In particular, the sequence $(x_n)$ is interpolating for $H^{\infty}(B_H)$ and the constant of interpolation is bounded by
$$\frac{128 (1+ 2 \log \frac{1}{\delta})}{e \delta }.$$
\end{teo}

\dem Define, for any $k,j \in \N$, $k \neq j$, the analytic function $g_{k,j}: H \To \C$ given by $g_{k,j}(x) = \langle \phi_{x_k}(x),\phi_{x_k}(x_j)\rangle.$ For each $j \in \N$ we define the function $B_j: B_H \To \C$ by $B_{j}(x) = \prod_{k \neq j} g_{k,j}(x).$ First we check that the infinite product converges
uniformly on $rB_{H}=\{x\in B_H:\|x\|\leq~r\}$ for fixed $0< r < 1$. Let $x \in r B_{H}$. We have, by Lemma \ref{lemacarlesonhilbert}, that
$$1-g_{k,j}(x)= 1-\langle \phi_{x_k}(x),\phi_{x_k}(x_j)\rangle=\frac{(1- \langle x,x_j \rangle)(1- \langle x_k,x_k \rangle)}{(1- \langle x,x_k \rangle)(1- \langle x_k,x_j \rangle)}.$$

It is easy that $|1- \la x,x_j \ra | \leq 1+r$, \ $|1- \la x,x_k \ra| \geq 1-r$ and $|1- \la x_k,x_j \ra| \geq  1-\|x_k\| \|x_j\| \geq 1-\|x_j\|.$ Then, we have that
$$|1-g_{k,j}(x)| \leq \frac{1+r}{1-r} \frac{1-\|x_k\|^2}{1-\|x_j\|},$$

\noindent so for any $j \in \N$, the series $\sum_{k \neq j} |1-g_{k,j}(x)|$ is uniformly convergent on $r B_H$ by Lemma \ref{lema4Hilbert}. In particular, the infinite product $\prod_{k \neq j} g_{k,j}(x)$ converges uniformly on compact sets, so $B_j \in H(B_H)$. In addition, notice that for $x \in B_H$, $|B_{j}(x)| = \prod_{k \neq j} |g_{k,j}(x)| = \prod_{k \neq j} |\langle \phi_{x_k}(x), \phi_{x_k}(x_j) \rangle| \leq \prod_{k \neq j} \|\phi_{x_k} (x)\| \| \phi_{x_k} (x_j)\| \leq 1$, so $\|B_{j}\|_{\infty} \leq 1$ and we obtain that $B_j \in H^{\infty}(B_H)$. It is clear that $B_{j}(x_{k}) =0$  for $k
\neq j$ since $\phi_{x_k}(x_k)=0$ and, according to \ref{metrica hilb}, we have that
$$|B_{j}(x_{j})| = \prod_{k \neq j} |g_{k,j}(x_j)| = \prod_{k \neq j} |\la \phi_{x_k}(x_j), \phi_{x_k}(x_j) \ra| =$$
$$\prod_{k \neq j} \| \phi_{x_k} (x_j)\|^2 = \prod_{k \neq j} \rho(x_k,x_j)^2 \geq \delta^2.$$

Consider the functions $q_j, A_j \in H(B_H)$ for any $k \in \N$ defined by
$$q_j(x)= \left(\frac{1-\|x_j\|^2}{1-\la x,x_j \ra}\right)^2,$$
%$$A_j(x)= \sum_{\left\{ k : \|x_{k}\| \geq \|x_{j}\| \right\}} (1- \rho(x_k,x_j)^2) \frac{1+<x,x_k>}{1-<x,x_k>}.$$
$$A_j(x)= \sum_{\left\{ k : \|x_{k}\| \geq \|x_{j}\| \right\}} \frac{(1-\|x_k\|^2)(1-\|x_j\|^2)}{1-|\la x_k,x_j\ra|^2} \frac{1+\la x_k,x\ra}{1-\la x_k,x\ra} .$$

The function $q_j$ is clearly analytic and bounded. By Lemma \ref{reordenar}, we will consider that the sequence $(\|x_n\|)$ is non-decreasing, so $\left\{ k : \|x_{k}\| \geq \|x_{j}\| \right\}=\{k : k \geq j\}$. Notice also that for $0 < r <1$ and $x \in r B_H$ we have that $|A_j(x)| \leq$
$$\sum_{k \geq j} \frac{(1-\|x_k\|^2)(1-\|x_j\|^2)}{1-\|x_j\|^2} \frac{1+r}{1-r} \leq \frac{1+r}{1-r} \sum_{k \geq j} (1-\|x_k\|^2)$$

\noindent so by Lemma \ref{lema4Hilbert}, the series converges uniformly on \ $r B_H$ \ and hence \ $A_j \in H(B_H)$. \ Moreover, $\exp \left(-A_j \right)$ belongs to $H^{\infty}(B_H)$ since $| \exp \left( - A_j \right)|= \exp \left( -\Re e \ A_j \right)$ and using formula \ref{lemametps3}, we have 
$$\Re e \ A_j(x) = \sum_{k \geq j} \frac{(1-\|x_k\|^2)(1-\|x_j\|^2)(1-|\la x_k,x\ra|^2)}{(1-|\la x_k,x_j\ra|^2)(|1-\la x_k,x\ra|^2)}>0.$$

Consider $C_\delta=1/(1+2 \log{1/\delta})$ and for any $j \in \N$, the analytic function $F_j: B_H \To \C$ defined by
$$F_j(x)= \frac{B_j(x)}{B_j(x_j)} q_j(x)^2 \exp{\left( -C_\delta(A_j(x)-A_j(x_j)) \right)}.$$

\noindent It is clear that $F_j(x_j)=1$ and $F_j(x_k)=0$ for any $k \neq j$. We claim that there exists $M>0$ such that
$\sum_{j=1}^{\infty} |F_j(x)| \leq M$ for any $x \in B_H$. Indeed, by (\ref{lemametps3}) that
$$\Re e \ A_j(x) = \sum_{k \geq j} \frac{(1-\|x_k\|^2)(1-\|x_j\|^2)(1-|\la x_k,x\ra|^2)}{(1-|\la x_k,x_j\ra|^2)(|1-\la x_k,x\ra|^2)}.$$

\noindent In particular, for $x=x_j$, we obtain
$$\Re e \ A_j(x_j) =\sum_{k \geq j} \frac{(1-\|x_k\|^2)(1-\|x_j\|^2)(1-|\la x_k,x_j\ra|^2)}{(1-|\la x_k,x_j\ra|^2)(|1-\la x_k,x_j\ra|^2)}.$$
%$$\sum_{k \geq j}  \right( 1- \rho(x_k,x_j)^2 \left).$$

Using formula (\ref{pseudometric Hilbert}) we obtain that
$$\Re e \ A_j(x_j) = \sum_{k \geq j} (1-\rho^2(x_k,x_j)) = 1+ \sum_{k > j} (1-\rho^2(x_k,x_j))$$

\noindent and by (\ref{lemametps2}), we have that $\Re e \ A_j(x_j) \leq$
\begin{eqnarray*} 
1- \sum_{\left\{ k : \|x_{k}\| > \|x_{j}\| \right\}} \log{\rho(x_k,x_j)^2} \leq 1- \sum_{k \neq j} \log{\rho(x_k,x_j)^2} \leq  1+ 2 \log{\frac{1}{\delta}}.
\end{eqnarray*}

%Moreover, since $1-|<x,x_k>|^2 \geq 1-\|x_k\|^2$ and $1-\|x_j\|^2 \geq 1-\|x_k\|^2$ for $k \in \N$ such that $\|x_k\| \geq \|x_j\|$, we obtain, making some calculations, that
%$$-\Re e \ (A_j(x)) \leq - \sum_{\left\{ k : \|x_{k}\| > \|x_{j}\| \right\}} \frac{(1-\|x_k\|^2)^3}{|1-<x,x_k>|^2}.$$

Moreover, to estimate $\Re e \ A_j(x)$ from below we use Lemma \ref{lema ber hilbert} and we obtain that
$$\Re e \ A_j(x) \geq \frac{1}{8} \frac{1-\|x_j\|^2}{1-|\la x_j,x\ra|^2} \sum_{k \geq j} \frac{(1-\|x_k\|^2)^2}{|1-\la x_k,x\ra|^2}.$$

We define for any $k \in \N$,
$$b_k(x)=\frac{1-\|x_k\|^2}{1-|\la x_k,x\ra|^2}$$

\noindent so
\begin{eqnarray} \label{parte real Aj2}
\Re e \ A_j(x) \geq \frac{1}{8} \  b_j(x) \sum_{k \geq j} |q_k(x)|.
\end{eqnarray}

It is clear that $1-|\la x_j,x \ra|^2= (1+|\la x_j,x \ra|)(1-|\la x_j,x \ra|) \leq 2 |1- \la x_j,x \ra|,$ so 
$$|q_j(x)| = \left| \frac{1-\|x_j\|^2}{1- \la x,x_j\ra} \right|^2 \leq 4 \left( \frac{1-\|x_j\|^2}{1-|\la x,x_j\ra|^2} \right)^2 = 4 b_j(x)^2.$$

Using that $|B_j(x_j)| \geq \delta$, $|B_j(x)| \leq 1$, the bound for $\Re e \ A_j(x_j)$ and \ref{parte real Aj2}, we obtain 
\begin{eqnarray*}
|F_j(x)| \leq \frac{4}{\delta} |q_j(x)| b_j(x)^2  \exp{\left( - C_{\delta} \left( \frac{1}{8} b_j(x) \sum_{k \geq j} |q_k(x)| -\frac{1}{C_\delta} \right) \right)} \leq \\
\frac{4e}{\delta} |q_j(x)|  b_j(x)^2  \exp{ \left(- \frac{1}{8} C_{\delta} b_j(x) \sum_{k \geq j} |q_k(x)| \right)}.
\end{eqnarray*}

\noindent Since $0 \leq b_k(x) \leq 1$, we consider $u=b_j(x)$ and $t=C_{\delta} \sum_{k \geq j} |q_k(x)|  > 0$ and apply Lemma \ref{lema ber hilbert 2} to conclude that 
$$|F_j(x)| \leq \frac{4e}{\delta C_\delta} C_\delta |q_j(x)| h \left( C_{\delta} \sum_{k \geq j} |q_k(x)| \right),$$

\noindent where $h(t)=\min\{1,256/e^2 t^2\}$. Hence, summing on $j$, we obtain
$$\sum_{j=1}^{\infty} |F_j(x)| \leq \frac{4e}{ \delta C_\delta} \sum_{j=1}^{\infty}  C_\delta |q_j(x)| h\left( \sum_{k \geq j} C_\delta |q_k(x)| \right),$$

\noindent  and applying Lemma \ref{func dec} , we obtain that 
$$\sum_{j=1}^{\infty} |F_j(x)| \leq \frac{4e}{\delta C_\delta} \int_{0}^\infty h(t) dt = \frac{4e.32}{e^2 \delta C_{\delta}}=\frac{128}{e \delta C_{\delta}}.$$

\noindent Hence , by Proposition \ref{prop 21}, we conclude that $(x_n)$ is linear interpolating. \qed \bigskip

Given $(x_n) \subset B_H$ satisfying the extended Carleson's condition and any $(\alpha_n) \in \ell_{\infty}$, the function $f(x)=\sum_{j=1}^\infty \alpha_j F_j(x)$, where $F_j$ is defined as in Theorem \ref{teo carl hilbert}, is well-defined and interpolates the values $\alpha_n$ in the points $x_n$ for any $n \in \N$.

Notice also that the function $$f(\delta)=\frac{1}{\delta C_\delta^2}=\frac{1+2 \log{1/\delta}}{\delta}$$

\noindent is non-increasing for $0 < \delta \leq 1$. Since $\lim_{\delta \rightarrow 1} f(\delta) =1$, un upper bound for the constant of interpolation is close to $\frac{128}{e}$ if we deal with sequences satisfying the extended Carleson's condition with $\delta$ close to $1$. Can the number $\frac{128}{e}$ be decreased?

%The \ sequence \ of \ Beurling \ functions in \ Theorem \ref{teo carl hilbert} \ is \ bounded \ by $\frac{2048}{e \delta} (1+ 2 \log{\frac{1}{\delta}})^2$, so this is an upper bound for the constant of interpolation of $(x_n)$.

\end{document}